# Ham Sandwich with Mayo:

# A Stronger Conclusion to the Classical Ham Sandwich Theorem

By John H. Elton and Theodore P. Hill

School of Mathematics, Georgia Institute of Technology

Atlanta GA 30332-0160 USA

elton@math.gatech.edu, hill@math.gatech.edu

**Summary**

The conclusion of the classical ham sandwich theorem of Banach and Steinhaus may be strengthened: there always exists a common bisecting hyperplane that touches each of the sets, that is, intersects the closure of each set. Hence, if the knife is smeared with mayonnaise, a cut can always be made so that it will not only simultaneously bisect each of the ingredients, but it will also spread mayonnaise on each. A discrete analog of this theorem says that $n$ finite nonempty sets in $n$-dimensional Euclidean space can always be simultaneously bisected by a single hyperplane that contains at least one point in each set. More generally, for $n$ compactly-supported positive finite Borel measures in Euclidean $n$-space, there is always a hyperplane that bisects each of the measures and intersects the support of each measure.

## 1. Introduction.

The classical ham sandwich theorem [BZ, S, ST] says that every collection of $n$ bounded Borel sets in $\mathbb{R}^n$ can be simultaneously bisected in Lebesgue measure by a single hyperplane. Many generalizations of this theorem are well known (e.g., [H], [M], [ST]), and the purpose of this note is to show that the conclusion of the classical ham sandwich theorem (and the conclusions of some of its well known extensions and generalizations) may be strengthened, without additional hypotheses. In the classical setting of bounded Borel sets, for example, it is shown that there always exists a Banach-Steinhaus bisecting hyperplane that contains at least one point in the closure of each of the sets. In the discrete setting, where the sets are finite, there always exists a bisecting hyperplane that contains at least one point in each of the sets. For compactly-supported positive finite Borel measures, there is always a hyperplane that bisects each of the measures and intersects the support of each measure. Note that to be able to treat these more general cases where it is possible that a hyperplane has positive measure, bisection of a measure has been defined in this paper (see below) to mean that no more than half the mass of the measure lies on either side of the hyperplane (not including the hyperplane); this is equivalent to the hyperplane being a median for the measures (see [H]).

*Remark*: The proof of Theorem 5 below for general finite measures only assumes the existence of a bisecting hyperplane for the case of purely atomic measures with finitely many atoms, so it also gives a proof of the existence of bisecting hyperplanes in the general case with the Borsuk-Ulam theorem having only been used for the purely atomic fintely-many atom case.

*Notation.* Fix $n \in \mathbb{N}$, and for $x, y \in \mathbb{R}^n$, let $|x|$ denote the Euclidean norm of *x*. For subsets *A* and *B* of $\mathbb{R}^n$, let $d(A, B)$ denote the Euclidean distance between *A* and *B*, i.e., $d(A, B) = \inf\{|x - y| : x \in A, y \in B\}$. Recall that every hyperplane *H* in $\mathbb{R}^n$ may be represented by $(u, c) \in \mathbb{R}^{n+1}$ via the relationship $x \in H \Leftrightarrow \langle u, x \rangle = c$, where *u* is a point in the unit *n*-sphere $\mathbb{S}^n$, $c \geq 0$ and $\langle \cdot, \cdot \rangle$ denotes the standard inner product on $\mathbb{R}^n$. For the hyperplane *H* determined by $(u, c)$, let $H^+$ denote the open half-space defined by $x \in H^+ \Leftrightarrow \langle u, x \rangle > c$ and $H^-$ denote the open half-space $x \in H^- \Leftrightarrow \langle u, x \rangle < c$.

For a bounded Borel set $A \subset \mathbb{R}^n$, let $\#\{A\}$ denote the cardinality of *A*, $\lambda(A)$ the Lebesgue measure (*n*-dimensional volume) of *A,* and $\overline{A}$ the closure of *A*. For a finite Borel measure $\mu$ on $\mathbb{R}^n$, $\|\mu\| = \mu(\mathbb{R}^n)$ denotes the total mass of $\mu$, and $supp(\mu)$ the support of $\mu$ (the smallest closed set $C \subset \mathbb{R}^n$ such that $\mu(C) = \|\mu\|$ ).

**Definition .** An *(n-1)*-dimensional *hyperplane H* in $\mathbb{R}^n$ *bisects a finite set* $A \subset \mathbb{R}^n$ if both $\#\{A \cap H^+\} \leq \#\{A\}/2$ and $\#\{A \cap H^-\} \leq \#\{A\}/2$; *bisects a bounded Borel set* $A \subset \mathbb{R}^n$ if $\lambda(A \cap H^+) = \lambda(A \cap H^-) = \lambda(A)/2$; and *bisects a positive finite Borel measure* $\mu$ on $\mathbb{R}^n$ if $\mu(H^+) \leq \|\mu\|/2$ and $\mu(H^-) \leq \|\mu\|/2$.

## 2. Bisecting Discrete Measures.

**Theorem 1**. *Let* $\mu_1, \ldots \mu_n$ *be purely atomic finite positive measures on* $\mathbb{R}^n$, *with finitely many atoms. Then there exists a hyperplane H such that for all i, H bisects* $\mu_i$ *and* $\mu_i(H) > 0$.

Although the bisection conclusion of Theorem 1 can be proved by first principles using the Borsuk-Ulam Theorem, the next lemma, a discrete version of the ham sandwich theorem, will facilitate its proof. The lemma follows easily from the classical ham sandwich theorem, and is a direct corollary of [H, Theorem 1].

**Lemma 2.** *Let $\mu_1,...\mu_n$ be purely atomic positive measures on $\mathbb{R}^n$ with finitely many atoms. Then there exists a hyperplane $H$ that bisects $\mu_i$ for all $1 \leq i \leq n$.*

*Proof of Theorem 1.* Fix $\varepsilon > 0$, and let $\{A_i\}$ be the sets of atoms of $\{\mu_i\}$, respectively. For each $i$, $1 \leq i \leq n$, reduce the mass of one of the atoms of $\mu_i$ by some positive amount less than $\varepsilon$ (and less than the mass of the smallest atom of $\mu_i$) such that for the new measure $\mu_i'$,

$$(*) \qquad \sum_{x \in S} \mu_i'(x) \neq \sum_{x \in A_i \setminus S} \mu_i'(x) \quad \text{for every } S \subset A_i.$$

(This can clearly be done since each $A_i$ is finite).

Now apply Lemma 2 to $\mu_1',...\mu_n'$. The resulting hyperplane $H_\varepsilon$ bisects each $\mu_i'$ and must in fact contain an atom of each $\mu_i'$ (which has the same atoms as $\mu_i$, just of different mass), or it could not bisect it, because of $(*)$. Let $\varepsilon$ approach zero along some sequence $\varepsilon_k$ such that the corresponding hyperplanes $H_k$ converge to say $H$. Clearly $H$ bisects each $\mu_i$, and since each $H_k$ contains an element of $A_i$ for each $i$ and the $A_i$ are finite, by passing to a subsequence it may be assumed that there is an atom of $\mu_i$ for each $i = 1,...n$ which belongs to all the $H_k$, and hence to $H$. Thus $\mu_i(H) > 0$, for each $1 \leq i \leq n$. □

**Corollary 3.** *For every collection $A_1,...,A_n$ of non-empty finite subsets of $\mathbb{R}^n$, there is a hyperplane $H$ such that for each $1 \leq i \leq n$, $H$ bisects $A_i$ and $H \cap A_i \neq \emptyset$.*

*Proof.* Let $\{\mu_i\}$ be the measures with atoms $\{A_i\}$, respectively, and masses of each atom equal to 1. Apply Theorem 1. □

**Example 4.** Sprinkle some salt and pepper on a table, any amounts of each. Then there is always a grain of salt and a grain of pepper and a line through both grains that has at most half of the grains of salt on each side, and also at most half of the grains of pepper on each side.

## 2. Bisecting General Measures.

**Theorem 5.** *For every collection $\mu_1,...\mu_n$ of compactly-supported positive Borel measures on $\mathbb{R}^n$ there exists a hyperplane $H$ such that for each $1 \leq i \leq n$, $H$ bisects $\mu_i$ and $H \cap \text{supp}(\mu_i) \neq \emptyset$.*

*Proof*: Let $C$ be a finite closed cube containing $\operatorname{supp}(\mu_i)$ for all i, and fix $\varepsilon > 0$. Let $P$ be a partition of $C$ into cubes (not necessarily closed or open) of diameter less than $\varepsilon$, and let $x_c$ be the centroid of cube $c$. For each i, let $v_i$ be the purely atomic measure such that for $c \in P, v_i(x_c) = \mu_i(c)$, and the only atoms are the $\{x_c\}$. That is, approximate the measures with purely atomic ones by concentrating all the mass at the centroids of the cubes, for those cubes in the partition which have non-zero mass. By Theorem 1, there is a hyperplane $H = H_\varepsilon$ such that for all $i$, $H$ bisects $v_i$, and $x_c \in H$ for some $x_c$ with $v_i(x_c) > 0$, so some point of support of $\mu_i$ lies within distance $\varepsilon$ of $H$; that is, $d(H, \operatorname{supp}(\mu_i)) < \varepsilon$. Let $A^+ = \cup\{c \in P : c \subset H^+\}$, so that $A^+$ is the union of the cubes of the partition that are entirely contained in $H^+$. Note that $\mu_i(A^+) = v_i(A^+)$, $\|\mu_i\| = \|v_i\|$, and since $H$ bisects $v_i$, $\mu_i(A^+) \leq v_i(H^+) \leq \|v_i\|/2 = \|\mu_i\|/2$. Note that any point in $H^+ \cap C$ whose distance from $H$ is greater than or equal to $\varepsilon$ belongs to $A^+$. $A^-$ is defined similarly; and similarly, $\mu_i(A^-) \leq \|\mu_i\|/2$.

Now let $\varepsilon = 1/k$, $k = 1, 2, 3, \ldots$, and let $H_k$, $A_k^+$ and $A_k^-$ correspond to $H$, $A^+$ and $A^-$ above. Since $C$ is compact, by passing to a subsequence if necessary, it may be assumed that the hyperplanes $H_k$ converge a hyperplane $H$, in such a way that $d(H_k \cap C, H) < 1/k$, and $(u_k, c_k) \to (u, c)$ where $(u_k, c_k) \in \mathbb{R}^{n+1}$ and $(u, c) \in \mathbb{R}^{n+1}$ represent $H_k$ and $H$, respectively, as in the earlier definition of hyperplanes. Also, $d(H, \operatorname{supp}(\mu_i)) < 1/k$, $\mu_i(A_k^+) \leq \|\mu_i\|/2$ and $\mu_i(A_k^-) \leq \|\mu_i\|/2$ for all $i$.

Note that $H^+ \subset (H^+ \setminus A_k^+) \cup A_k^+$, so $\mu_i(H^+) \leq \|\mu_i\|/2 + \mu_i(H^+ \setminus A_k^+)$ for all $k$. It will be shown below that the sets $(H^+ \cap C) \setminus A_k^+ \to \emptyset$, so from the continuity theorem for measures, $\mu_i(H^+ \setminus A_k^+) = \mu_i((H^+ \cap C) \setminus A_k^+) \to 0$, and therefore $\mu_i(H^+) \leq \|\mu_i\|/2$; and from $d(H_k, \operatorname{supp}(\mu_i)) < 1/k$ it follows that $H \cap \operatorname{supp}(\mu_i) \neq \emptyset$ since $\operatorname{supp}(\mu_i)$ is closed. Similarly $\mu_i(H^-) \leq \|\mu_i\|/2$. This will finish the proof, once it is shown that $(H^+ \cap C) \setminus A_k^+ \to \emptyset$.

Now $H^+ \setminus H_k^+ \to \emptyset$, because $\langle u, x \rangle > c \Rightarrow \langle u_k, x \rangle > c_k$ for sufficiently large $k$, so $x \in H^+ \Rightarrow x \in H_k^+$ for sufficiently large $k$. Suppose $x \in (H_k^+ \cap C \cap H^+) \setminus A_k$. Then since $x \notin A_k$, $d(x, H_k) < 1/k$ from the definition of $A_k$. Since $d(H_k \cap C, H) < 1/k$, it follows that $d(x, H) < 2/k$. So if $x \in (H_k^+ \cap C \cap H^+) \setminus A_k$ for all $k$, it would follow that $x \in H$, which is impossible since $H$ is disjoint from $H^+$. So $(H_k^+ \cap C \cap H^+) \setminus A_k \to \emptyset$ also.

Thus since $(H^+ \cap C) \setminus A_k^+ \subset (H^+ \setminus H_k^+) \cup (H_k^+ \cap C \cap H^+) \setminus A_k$, $(H^+ \cap C) \setminus A_k^+ \to \emptyset$ as claimed. $\square$

**Example 6.** At any given instant of time, there is one planet, one moon and one asteroid in our solar system and a single plane touching all three that exactly bisects the total planetary mass, the total lunar mass, and the total asteroidal mass of the solar system. (Note that different objects may have different mass densities, and even non-uniform mass densities, so this conclusion does not follow from the next corollary.)

**Corollary 7. (Ham Sandwich with Mayo).** *For every collection $A_1,..., A_n$ of n bounded Borel subsets of $\mathbb{R}^n$ of positive Lebesgue measure, there exists a hyperplane H such that for each $1 \leq i \leq n$, H bisects $A_i$ and $H \cap \overline{A_i} \neq \emptyset$.*

*Proof.* The conclusion follows immediately from Theorem 5 by letting $\mu_1,...,\mu_n$ be the finite Borel measures on $\mathbb{R}^n$ defined by $\mu_i(B) = \lambda(B \cap A_i)$ for all Borel sets $B \subset \mathbb{R}^n$, and observing that $\text{supp}(\mu_i) \subset \overline{A_i}$. □

If the sets $A_1,..., A_n$ are all closed, of course, then there is always a bisecting hyperplane that intersects each set. Otherwise, as the following simple example shows, there may not be a bisecting hyperplane that intersects any of the sets.

**Example 8.** Let $n = 2$, and suppose that $A$ is the union of the two open disks of radius one centered at (-1,1) and at (-1,-1), and $B$ is the union of the two open disks of radius one centered at (1,1) and at (1,-1). It is easy to see that the unique line that bisects both $A$ and $B$ is the line $y = 0$, which intersects the closure of both $A$ and $B$, but does not intersect either set.

Even if the sets are closed, not all of the bisecting hyperplanes guaranteed to exist by the classical ham sandwich theorem will intersect all the bisected sets.

**Example 9.** Let $n = 2$, and suppose that $A$ is the closed disk of radius one centered at (0,0), and $B$ is the union of the two closed disks of radius one centered at (-3,0) and at (3,0). Then the line $x = 0$ bisects both $A$ and $B$, but does not touch $B$. The line $y = 0$ bisects and intersects both sets.

If the hypothesis of compactly-supported is dropped in Theorem 5, the bisection conclusion still holds (cf. [H]), but there may not be a common bisecting hyperplane that intersects the supports of all the measures.

**Example 10.** Let $n = 2$, $A_1 = \{(k,1) \in \mathbb{R}^2 : k \in \mathbb{N}\} \cup \{(k,-1) \in \mathbb{R}^2 : k \in \mathbb{N}\}$, $A_2 = \{(-k,2) \in \mathbb{R}^2 : k \in \mathbb{N}\} \cup \{(-k,-2) \in \mathbb{R}^2 : k \in \mathbb{N}\}$, and let $\mu_1$ and $\mu_2$ be the purely atomic probability measures with supports $A_1$ and $A_2$, respectively, given by $\mu_1(k,1) = \mu_1(k,-1) = \mu_2(-k,2) = \mu_2(-k,-2) = 2^{-(k+1)}$ for all $k \in \mathbb{N}$. It is easy to see that the only lines that bisect both measures simultaneously are horizontal lines with height in [-1,1]. But none of these intersects $A_2$, the support of $\mu_2$.

**Correction.** The proof of the *General Case* of Theorem 1 in [H] has a gap, namely, that although the subsequence of $(a_j, b_j) \in \mathbb{R}^{n+1}$ converges to $(a,b) \in \mathbb{R}^{n+1}$, it could be that $a = 0$ or $b = \infty$, in which case $\{x \in \mathbb{R}^n : \langle a, x \rangle = b\}$ does not define a hyperplane. To complete the argument, first note that no $a_j$ can be zero, so without loss of generality the $a_j$'s all have norm 1 (this may violate the normalization of the vectors $(a_j, b_j)$, but that is not a problem); thus $\|a\| = 1$. Since the $\{\mu_i\}$ are all probability measures, there exists an $r^*$ so that $\mu_1(\{x \in \mathbb{R}^n : \|x\| \geq r^*\}) < 1/2$. Replacing $(a_j, b_j)$ by $(-a_j, -b_j)$ where necessary, assume without loss of generality that all the $b_j$ are nonnegative. By construction, $\mu_1(\{x \in \mathbb{R}^n : \langle x, a_j \rangle \geq b_j\}) \geq 1/2$ for all *j*, so by Schwarz's inequality, $\mu_1(\{x \in \mathbb{R}^n : \|x\| \geq b_j\}) \geq 1/2$. Hence by the definition of $r^*$, $b_j \leq r^*$ for all *j*, so $b < \infty$. □

**Acknowledgment.** The second author is grateful to Professor Robert Burckel for pointing out the gap noted above, and for suggesting arguments to correct it.

### References.